\documentclass[12pt,a4paper]{article}
\usepackage{amsmath,amsfonts,amssymb,amsthm,amscd}
\usepackage{hyperref}
\usepackage{epstopdf}

\usepackage{braket}
\usepackage{xypic}
\usepackage{amsmath}

\usepackage[american]{babel}
\usepackage{color} 
\usepackage{epsfig}
\usepackage{caption}
\usepackage{rotating}
\usepackage{setspace}
\usepackage{fancyhdr}
\usepackage{booktabs}
\usepackage{mathrsfs}
\usepackage{amsthm}
\usepackage{amsmath,amssymb}
\usepackage{enumerate}
\usepackage{tikz}
\usepackage{youngtab}
\usepackage{hyperref}
\usepackage{tensor}
\usepackage{mathabx}
\usepackage{booktabs}
\usepackage{mathrsfs}
\usepackage{subfigure}
\usepackage{amsmath,amssymb}

\renewcommand{\[}{\begin{equation}}
\renewcommand{\]}{\end{equation}}

\newtheorem{thm}{Theorem}[section]
\newtheorem{cor}[thm]{Corollary}
\newtheorem{lem}[thm]{Lemma}
\newtheorem{prop}[thm]{Proposition}

\theoremstyle{definition}
\newtheorem{defn}[thm]{Definition}

\theoremstyle{remark}
\newtheorem{oss}[thm]{Remark}

\newcommand{\ra}{\rightarrow}

\newcommand{\nn}{\mathbb{N}}
\newcommand{\cc}{\mathbb{C}}

\DeclareMathOperator{\id}{id}

\DeclareMathOperator{\tr}{tr}

\DeclareMathOperator{\Aut}{Aut}

\DeclareMathOperator{\diag}{diag}

\DeclareMathOperator{\ad}{ad}
\DeclareMathOperator{\Iso}{Iso}
\DeclareMathOperator{\ISO}{ISO}
\DeclareMathOperator{\chara}{Char}

\DeclareMathOperator{\dom}{dom}

\title{On isometries of spectral triples associated to $AF$-algebras and crossed products }

\begin{document}

\author{Jacopo Bassi and Roberto Conti}
\date{}
\maketitle



\begin{abstract}
\noindent We examine the structure of two possible candidates of isometry groups for the spectral triples on $AF$-algebras introduced by Christensen and Ivan. In particular, we completely determine the isometry group introduced by Park, and observe that these groups coincide in the case of the Cantor set. We also show that the construction of spectral triples on crossed products given by Hawkins, Skalski, White and Zacharias, is suitable for the purpose of lifting isometries.




\end{abstract}
\section{Introduction}
The reconstruction theorem by Connes is a smooth version of the Gelfand duality; it shows that both the geodesic distance and the volume form on a compact Riemannian space are encoded in a particular spectral triple (\cite{Co}). The non-commutative differential geometry paradigm is that spectral triples on $C^*$-algebras should allow the investigation of geometric properties in the non-abelian setting.\\ Many authors studied the nature of such features in different specific cases; this lead for example to the notion of compact quantum metric space, as defined by Rieffel (\cite{Ri}). Examples of non-commutative manifolds have been studied in the case of the reduced $C^*$-algebra of a discrete group (\cite{Co}), $AF$-algebras (\cite{AnCh2}), crossed products (\cite{HSWZ}), and many other instances.\\
From the point of view of non-commutative geometry, a natural problem is to understand what is the right notion of isometry for a non-commutative manifold. Such topic was investigated in \cite{Pa1, Pa2, Pa3, LoWu, CoRo} and at the time being the two definitions of non-commutative isometry appearing in these manuscripts are the only natural candidates known to the authors. 
More precisely, these are the automorphisms implemented by unitaries commuting with the Dirac operator and the automorphisms leading to preservation of the Connes distance on the state space, giving rise to the two groups $\Iso$ and $\ISO$, respectively.


In the present work the authors investigate properties of the $\Iso$- and $\ISO$-group in some specific cases, namely the spectral triples constructed by Christensen and Ivan on $AF$-algebras and the spectral triples 
on the crossed product of a $C^*$-algebra with a discrete group, as defined in \cite{HSWZ}. A concrete description of the $\Iso$-group is provided in the case of the spectral triples for $AF$-algebras (under a natural nondegeneracy condition), namely its elements are the automorphisms which preserve the filtration and the given reference faithful state. In the particular case of some natural spectral triple on the Cantor set, this group actually coincides with the $\ISO$-group.
Moreover, it is proved that in the case of $UHF$-algebras of type $n^\infty$, if the eigenvalues of the Dirac operator grow fast enough, the $\ISO$-group cannot contain all the switches of the tensor factors. In the case of the CAR-algebra, we explicitly compute the Connes distance between certain states; as a consequence,
none of the switch automorphisms can appear in the $\ISO$-group at all. In Section 4 it is shown a procedure to lift the elements of the $\Iso$-group of a spectral triple for a $C^*$-algebra endowed with an action of a countable discrete group to the $\Iso$-group of the spectral triple of the crossed product given in \cite{HSWZ}.

\section{Preliminaries}
 
\subsection{The definition of a spectral triple}

The geometric properties of a compact Riemannian spin manifold can be reconstructed from the algebraic data contained in the way the Dirac operator and the measurable bounded functions interact when acting on the Hilbert space of $L^2$-spinors (\cite{Co}). These data are encoded in the notion of spectral triple, which we recall in the setting of arbitrary unital $C^*$-algebras.

\begin{defn}
\label{defspec}
Let $A$ be a unital $C^*$-algebra represented on a Hilbert space $H$ and $D$ be an unbounded self-adjoint operator on $H$. Denote $A_D= \{ a \in A \; | \; a \dom (D) \subset \dom (D)\}$. We say that $(A, H, D)$ is a {\it spectral triple} if the following conditions are satisfied:
\begin{itemize}
\item[(i)] the set $\{a \in A_D \; | \; [D,a] \mbox{ extends to a bounded operator on } H\}$ is norm-dense in $A$;
\item[(ii)] The operator $(1+D^2)^{-1}$ is compact.
\end{itemize}
In this case, the operator $D$ is a \it{Dirac operator}.
\end{defn}
As observed by Connes, given a spectral triple $(A,H,D)$, it is possible to define a pseudo-metric $d_D$ on $\mathcal{S}(A)$, the state space of $A$, by the formula 
\begin{equation*}
d_D (\phi, \psi)= \sup_{ a \in A \; : \; \|[D,a]\| \leq 1} |\phi (a) - \psi (a)|.
\end{equation*}
In the case where the spectral triple is the natural one associated to a compact Riemannian spin manifold, this formula gives back the geodesic distance (\cite{Co}, chapt.~6, par.~1); there are also examples of non-commutative $C^*$-algebras admitting spectral triples for which $d_D$ is a metric on the state space inducing the $w^*$-topology (c.f.~\cite{Ri}).

\subsection{Spectral triples on $AF$-algebras}
We briefly recall the construction of spectral triples on $AF$-algebras given in \cite{AnCh2}. Let $A= \lim_n A_n$ be a unital $AF$-algebra associated to a filtration of finite-dimensional $C^*$-algebras $A_n$ and assume $A_0= \mathbb{C} 1$. Let $\phi$ be a faithful state on $A$ and $H_\phi$ the associated $GNS$-Hilbert space, $\xi$ the corresponding cyclic vector. For every $n \in \mathbb{N}$ let $H_n = A_n \xi$ and $P_n : H_\phi \rightarrow H_n$ be the associated orthogonal projection; we define $K_n = H_n \ominus H_{n-1}$, $Q_n =P_n - P_{n-1}$. Given a sequence of positive real numbers $\lambda_n \rightarrow \infty$, with $\lambda_0 =0$, we let $D_{\{\lambda_n\}}$ be the unbounded self-adjoint operator on $H_\phi$ given by $D=\sum_n \lambda_n Q_n$ with its natural domain; this is a Dirac operator for a spectral triple $(A, H_\phi, D_{\{\lambda_n\}} )$ on $A$. The elements of the dense subalgebra $\bigcup_n A_n \subset A_{D_{\{\lambda_n\}}}$ have bounded commutators with this Dirac operator. Following \cite{HSWZ}, spectral triples constructed in this way will be referred to as Christensen-Ivan spectral triples and we will drop the suffix $\{\lambda_n\}$ appearing in the Dirac operator where no confusion is likely to arise. As shown in \cite{AnCh2} Theorem 2.1, if the sequence of eigenvalues $\{\lambda_n\}$ grows rapidly enough, the pseudo-metric induced by $D_{\{\lambda_n\}}$ on the state space of $A$ is actually a metric and it induces the $w^*$-topology. We will see an application of this fact in Lemma \ref{leml}.

\subsection{Isometries of non-commutative spaces}
Park introduced in \cite{Pa1} the concept of isometry in the non-commutative setting (see also \cite{LoWu}). He showed that in the case of the spectral triple given by the continuous functions on a compact Riemannian oriented manifold, the Hilbert space of complex $L^2$-forms and the de Rham operator as Dirac, this concept coincides with the ordinary notion of isometry. He also studied the group of such isometries in some non-commutative cases. 
\begin{defn}[\cite{Pa1}]
Let $(A, H, D)$ be a spectral triple. An element $\alpha \in \Aut (A)$ belongs to the group $\Iso (A,H,D)$ if $\alpha$ is implemented on $H$ by a unitary operator which leaves the domain of the Dirac operator $D$ invariant and commutes with $D$.
\end{defn} 
In view of the possibility to characterize the distance between two points in a compact Riemannian spin manifold $M$ in terms of the metric $d_D$ on the state space of $C(M)$, another possible non-commutative generalization of isometry which is worth studying is given by:
\begin{defn}[\cite{CoRo}]
Let $(A,H,D)$ be a spectral triple. An element $\alpha \in \Aut (A)$ belongs to the group $\ISO (A,H,D)$ if $d_D (\phi \circ \alpha, \psi \circ \alpha) = d_D (\phi, \psi)$ for every $\phi, \psi \in \mathcal{S}(A)$.
\end{defn}
While it is always the case that $\Iso (A,H,D) \subset \ISO (A,H,D)$, there are simple examples in which these two groups do not coincide.

\begin{prop}
\label{propflip}
Let $A=M_2(\mathbb{C})$ and let $D$ be a diagonal self-adjoint element of $A$ with two distinct eigenvalues. Then $\Iso (A,  \mathbb{C}^2, D) \neq \ISO(A,  \mathbb{C}^2, D)$.
\end{prop}
\proof The "flip" unitary $U =(1 - \delta_{i,j})_{i,j=1}^2$ does not commute with $D$, from which it follows that $\ad (U) \notin \Iso (A,\mathbb{C}^2, D)$; we will see that $U$ implements an element of $\ISO(A, \mathbb{C}^2, D)$. For let $\phi$, $\psi \in \mathcal{S}(A)$ and note that
\begin{equation*}
\begin{split}
d_D(\phi \circ \ad(U), \psi \circ \ad (U)) &= \sup_{ a \in A: \; \|[D,a]\| \leq 1} |\phi (UaU^*) - \psi (UaU^*)|\\
&= \sup_{a \in A: \; \|[D, U^* a U]\| \leq 1} |\phi (a) - \psi(a)|.
\end{split}
\end{equation*}
Hence it is enough to show that $\|[ D, UaU^*]\| = \|[D,a]\|$ for every $a \in A$. This follows from the equality $[D, U^* a U] = -U^*[D, a]U$ for every $a \in A$. \hfill $\Box$

\begin{oss} 
The fact that the automorphism induced by the flip unitary of Proposition \ref{propflip} preserves the Connes distance between \textit{pure} states already follows from the computation in \cite{IKM} Proposition 2.
\end{oss}

Other comparison results concerning the two isometry groups for the case of Cuntz algebras have been obtained in \cite{CoRo}.

\section{Isometries of Christensen-Ivan spectral triples}
This section contains information about the isometry groups associated to spectral triples on $AF$-algebras. We completely characterize the $\Iso$-groups for the cases at hand and collect some results concerning the "size" of the $\ISO$-group.

Let $A=\overline{\bigcup_n A_n}$ be an $AF$ algebra, $\phi$ a faithful state on $A$ and $(\pi_\phi , H_\phi, \xi_\phi)$ the associated GNS triple. From now on we identify $A$ with its image under $\pi_\phi$ and write $\xi$ for the cyclic and separating vector $\xi_\phi$. Let $(A,H_\phi,D )$ be the Christensen-Ivan spectral triple associated to the given generating family of finite-dimensional subalgebras $A_n$ and the diverging sequence of non-negative real numbers $\{\lambda_n\}$ (here, $A_0 = \cc I$ and $\lambda_0 = 0$).
\begin{thm}
\label{p1}
Let $(A,H_\phi , D)$ be as above and suppose that the sequence $\{\lambda_n\}$ satisfies $\lambda_i \neq \lambda_j$ for every $i \neq j$. Then
\begin{equation*}
\Iso(A,H_\phi,D) = \{\alpha \in \Aut(A) \ | \ \alpha(A_i) = A_i \ \forall i, \ \phi \circ \alpha = \phi\}
\end{equation*}
In particular, if $A= \otimes_{i=0}^\infty M_{n_i}$, with $M_{n_0} = \cc$, is a UHF-algebra, then
\begin{equation*}
\Iso (A,H_\phi,D)=\{ \alpha \in \Aut(A) \; | \; \alpha=\otimes_{i=1}^\infty \alpha_i, \; \alpha_i \in \Aut(M_{n_i}), \; \phi \circ \alpha = \phi\}
\end{equation*}
and if $\phi = \otimes_i \phi_i$ is a product state, then
\begin{align*}
\Iso (A,H_\phi,D) & = 
\{ \alpha \in \Aut(A) \; | \; \alpha=\otimes_{i=1}^\infty \alpha_i, \; \alpha_i \in \Aut(M_{n_i}), \; \phi_i \circ \alpha_i = \phi_i   \; \textit{ for every } i \in \nn\} \\
& \simeq \prod_{i=1}^\infty\{ u_i \in U_{n_i} \; | \; \phi_i \circ {\rm ad}u_i = \phi_i \; \textit{ for every } i \in \nn\}/ S^1.
\end{align*}
\end{thm}
\proof 
Let $\alpha \in \Iso (A,H_\phi,D)$ and $U \in U(H_\phi)$ be such that $\alpha= \ad (U)$. Since the coefficients $\lambda_n$ are pairwise different, the condition that $[D,U]=0$ is equivalent, for $U \in U(H)$, to the condition that $UH_i \subset H_i$ for every $i \in \nn$, where $H_i$ is the 
finite-dimensional Hilbert space $A_i \xi$, and thus $U H_i = H_i$, for every $i$. In particular, $U \xi \in \cc\xi$ and, without loss of generality, we may even assume that $U \xi = \xi$. Anyway, it readily follows that $\phi \circ \alpha = \phi$. Moreover,
note that, if $U \in B(H)$ leaves every $H_i$ invariant, then for every $i \in \nn$ and $a \in A_i$ we have $UaU^*\xi \in H_i$. Since $\xi$ is separating it follows that $UaU^* \in A_i$; hence $\ad (U)$ respects the specified filtration of $A$. Hence $\Iso (A,H_\phi,D) \subset \{\alpha \in \Aut(A) \ | \ \alpha(A_i) = A_i \ \forall i, \ \phi \circ \alpha = \phi\}$. Let now $\alpha$ be an automorphism of $A$ which leaves every $A_i$ invariant and preserves the faithful state $\phi$ under precomposition. Then the linear map $a\xi \ra \alpha(a)\xi$ defines a unitary operator on $H$; hence $\alpha$ is unitarily implemented and since the unitaries realizing $\alpha$ leave every $H_i$ invariant, they commute with $D$. 

Suppose now $A= \otimes_{i=0}^\infty M_{n_i}$ be a UHF-algebra and let $\alpha \in \Iso(A,H_\phi,D)$ be implemented by $U \in B(H_\phi)$. Since for every $i \in \nn$ $\ad (U)$ leaves both $M_{n_i}$ and $M_{n_i} \otimes M_{n_{i+1}}$ invariant, it also leaves $M_{n_{i+1}}$ invariant, being the commutant of $M_{n_i}$ in $M_{n_i} \otimes M_{n_{i+1}}$. Hence every $\alpha \in \Iso(A,H_\phi,D)$ decomposes as $\alpha = \otimes_{n=1}^\infty \alpha_i$, where $\alpha_i \in \Aut(M_{n_i})$. On the other hand, every automorphism of the form $\alpha = \otimes_{i=1}^\infty \alpha_i$ with $\alpha_i \in \Aut(M_{n_i})$ which satisfies $\phi \circ \alpha = \phi$ belongs to $\Iso(A,H_\phi,D)$. Now the last statement easily follows. \hfill $\Box$ 

The above argument can be used to obtain similar results in the more general case of spectral triples associated to non-decreasing sequences of real numbers.

\begin{cor}
\label{c1}
Let $A= \lim_i A_i$ be an AF-algebra, $\phi$ a faithful state on $A$ and $\lambda_i$ be a monotone non-decreasing sequence associated to a Dirac operator on $H_\phi$ implementing a Christensen-Ivan spectral triple. For every $i \in \nn$ let $n_i := \max \{k \geq i \; | \; \lambda_k = \lambda_i\}$ and $k_i := \min \{k \leq i \; | \; \lambda_k = \lambda_i\}$. Then
\begin{equation*}
\Iso (A,H_\phi,D) = \{ \alpha \in \Aut(A) \; | \; \alpha (A_{n_i}) = A_{n_i} \; \forall i \in \nn, \; \phi \circ \alpha = \phi  \}.
\end{equation*}
If $A=  \otimes_{i=0}^\infty M_{n_i}$ is UHF, then
\begin{equation*}
\Iso (A,H_\phi,D)=\{ \alpha \in \Aut(A) \; | \; \alpha=\otimes_{\lambda_i} \alpha_{\lambda_i}, \; \alpha_{\lambda_i} \in \Aut(\otimes_{j=k_i}^{n_i} M_{j}), \; \phi \circ \alpha_{\lambda_i} = \phi\}
\end{equation*}
\end{cor}
\proof
The Dirac operator reads $D= \sum_{\lambda_i \in \{\lambda_i\}_{i \in \nn}} \lambda_i(P_{n_i} - P_{k_i})$, where for $j \in \nn$, $P_j = P_{H_j}$. Hence $\Iso(A,H_\phi,D)$ can be computed using Theorem \ref{p1} applied to the $AF$-structure of $A$ given by the filtration $A=\lim_i A_{n_i}$. The second statement follows from the second part of the proof of Theorem \ref{p1}. \hfill $\Box$

\begin{oss}
In the particular case of a UHF-algebra $A= \otimes_{i=0}^{\infty} M_{n_i}$ represented on the GNS-Hilbert space associated to the  unique trace $\tau$, Corollary \ref{c1} gives, for a monotone non-decreasing sequence $\{\lambda_i\}$,
\begin{equation*}
\Iso(A,H_{\tau}, D_{\{\lambda_i\}})=\prod_{\lambda_i} U_{\prod_{\{j\; : \; \lambda_j=\lambda_i\}} n_j} /S^1
\end{equation*}
\end{oss}
 
\begin{oss}
It is unlikely that $\Iso (A,H_\phi,D)$ be normal in $\Aut (A)$. It is certainly not for $A= \otimes_{i=0}^{\infty} M_{n_i}$ with $n_i = m^k$ for some $2 \leq \min\{m,k\} $, for all $i \in \mathbb{N}$ and $\phi= \tr$. In fact, suppose $A=\otimes_{i=0}^\infty M_{2^k}$ and consider $\alpha= \id^{\otimes j} \otimes \hat{\alpha} \otimes \id^{\infty}$ for some $j \geq 1$, where $\hat{\alpha} (a_1 \otimes a_2 \otimes ... \otimes a_k) = a_2 \otimes a_1 \otimes a_3 \otimes ... \otimes a_k$. Then $\beta \circ \alpha \circ \beta^{-1} \notin \Iso (A, H_\phi , D)$ for $\beta = \id^{\otimes (j-1)} \otimes \hat{\beta} \otimes \id^{\otimes \infty}$, with $\hat{\beta} ( (a_1 \otimes ... \otimes a_k) \otimes (b_1 \otimes ... \otimes b_k)) = (a_1 \otimes ... \otimes a_{k-1} \otimes b_1 ) \otimes (a_k \otimes b_2 \otimes ... \otimes b_k)$.\end{oss} 

It was observed in \cite{AnCh2} that in the case of Christensen-Ivan spectral triples associated to rapidly diverging sequences of real numbers, all the information about the $w^*$-topology on the state space of the $AF$-algebra is carried by the Dirac operator, namely the topology induced by the metric associated to $D$ (see the discussion after Definition \ref{defspec}) on the state space of $A$ coincides with the $w^*$-topology (\cite{AnCh2} Theorem 2.1 (i)). The following is an application of this fact.

\begin{prop}
\label{p2}
Let $A$ be a UHF-algebra, $\tr$ the trace on $A$ and $\{\lambda_n\}_{n \in \mathbb{N}}$ a diverging sequence of positive real numbers such that the topology induced by the metric associated to $D_{\{\lambda_n\}}$ on $\mathcal{S}(A)$ coincides with the $w^*$-topology. Let $(A,H_{\tr},D_{\{\lambda_n\}})$ be the Christensen-Ivan spectral triple associated to this data. Then
\begin{equation*}
 \ISO(A,H_{\tr},D_{\{\lambda_n\}}) \subsetneq \Aut(A)
\end{equation*}
Moreover, $\Aut (A) \backslash \ISO (A, H_{\tr},D_{\{\lambda_n\}})$ contains $\ISO(A,H_{\tr},D_{\{\lambda_n\}}) \times \mathbb{N}$ as a subset.
\end{prop}
\proof In virtue of \cite{AnCh2} Theorem 2.1 (i) the pseudometric $d_{D_{\{\lambda_n\}}}$ induced by $D_{\{\lambda_n\}}$ on the state space of $A$ is actually a metric and it induces the $w^*$-topology. By \cite{Po} Corollary 3.8 the action of $\Aut (A)$ on the pure states of $A$ is transitive and by \cite{Gl} Theorem 2.8 the set of pure states is $w^*$-dense in the set of states. Hence for every pure state $\psi$ there is a sequence $\{\beta_n\} \subset \Aut (A)$ such that $d_{D_{\{\lambda_n\}}} (\psi \circ \beta_n , \tr) \rightarrow 0$. In particular, there is a sequence of positive real numbers $r_n \rightarrow 0$, $r_n \neq r_m$ for $n \neq m$ such that the sets $\Omega_n := \{ \phi \in \mathcal{S}(A) \; | \; d_{D_{\{\lambda_n\}}} (\tr , \phi) = r_n\} \neq \emptyset$ for every $n \in \mathbb{N}$ and from the above, given $n\neq m \in \mathbb{N}$, there is always an element $\alpha_{n,m} \in \Aut (A)$ such that $\alpha_{n,m} (\Omega_n) \cap \Omega_m \neq \emptyset$. We claim that given $n,m \in \mathbb{N}$ there is $l \in \mathbb{N}$ such that $\alpha_{n,m} \neq \alpha_{n,h} \circ \beta$ for every $h \geq l$, $\beta \in \ISO(A,H_{\tr}, D_{\{\lambda_n}\})$; indeed, if this is not the case, then it is possible to find a sequence $\phi_k \in \Omega_n$ such that $\alpha_{n,m} (\phi_k) \ra \tr$, by compactness we can suppose $\phi_k$ converges to a state $\phi$, which satisfies $\alpha_{n,m} (\phi)= \tr$, which is impossible. It follows that for every $n \in \mathbb{N}$ there is a sequence $m_i$ such that the elements $\alpha_{n,m_i} \circ \beta_1 \neq \alpha_{n,m_j} \circ \beta_2$ for every $(i,\beta_1)\neq (j,\beta_2) \in \mathbb{N} \times \ISO(A,H_{\tr}, D_{\{\lambda_n\}})$. \hfill $\Box$

Let $n \in \mathbb{N}$ and $A=\otimes^\infty M_n$ be a $UHF$-algebra. For $i,j \in \mathbb{N}$ let $\sigma_{i,j}$ be the automorphism of $A$ given by the permutation of the $i$-th and $j$-th tensor components. It follows from Theorem \ref{p1} that for every $i \neq j$, $\sigma_{i,j}$ does not belong to the $\Iso$-group of a Christensen-Ivan spectral triple for $A$ with pairwise different eigenvalues.

\begin{lem}
\label{leml}
Let $A= \otimes_i M_{n_i}$ be a UHF-algebra with $n_i =n_j=n$ for every $i,j$ and consider the setting as in the statement of Proposition \ref{p2}. There exists $i \in \mathbb{N}$ such that $\sigma_{i,i+1} \notin \ISO (A,H_{\tr},D_{\{\alpha_n\}})$.
\end{lem}
\proof
Let $a \in M_{n}$ be such that $\tr (a^* a)=1$, $\omega_a$ the associated (vector) state, namely $\omega_a (x) = \tr ((a^* \otimes 1)x(a \otimes 1))$ and suppose $\omega_a \neq \tr$; let $\phi : A \rightarrow A$ be the shift endomorphism: $\phi (x)= 1 \otimes x$. Then $\omega_{\phi^k a} \rightarrow \tr$ in the $w^*$-topology when $k \to \infty$. Since this topology coincides with the topology induced by the metric $d_{D_{\{\lambda_n\}}}$, we have $d_{D_{\{\lambda_n\}}} (\tr, \omega_{\phi^k a}) \rightarrow 0$, $k \to \infty$. Hence there is $k \in \mathbb{N}$ such that $d_{D_{\{\lambda_n\}}} (\tr, \omega_{\phi^k a}) < d_{D_{\{\lambda_n\}}} (\tr, \omega_a)$. Since $\omega_{\phi^k a} = \omega_a \circ \sigma_{1,k+1}$ and $\sigma_{1,k+1} = \sigma_{1,2} \circ \sigma_{2,3} \circ ...  \circ \sigma_{k-1,k} \circ \sigma_{k,k+1} \circ \sigma_{k-1,k} \circ ... \circ \sigma_{2,3} \circ \sigma_{1,2}$, the result follows. \hfill $\Box$\\


\begin{oss}
Let $A=\otimes^{\infty} M_n$ for some $n \in \nn$ and let $\{\lambda_i\}$ be a strictly increasing sequence of positive real numbers, $D_{\{\lambda_i\}}$ the associated Dirac operator. Then $\Iso (A, H_{\tr}, D_{\{\lambda_i\}})$ satisfies:
\begin{itemize}
\item[(i)] $\sigma_{i,j} \notin \Iso(A,H_{\tr}, D_{\{\lambda_i\}})$ for every $i,j$;
\item[(ii)] $\sigma_{i,j}$ is in the normalizer of $\Iso (A,H_{\tr}, D_{\{\lambda_i\}})$ in $\Aut (A)$ for every $i,j$;
\item[(iii)] 
 if $\alpha \in \Aut (A)$ then  $\alpha \in \Iso(A,H_{\tr}, D_{\{\lambda_i\}})$ if and only if $\id_{M_n}  \otimes \alpha \in \Iso(A,H_{\tr}, D_{\{\lambda_i\}})$.
\end{itemize}
\end{oss}

We will see below how a stronger version of Lemma \ref{leml} can be deduced in the concrete case of the CAR-algebra. Before specializing, we prove a result for the case of more general $UHF$-algebras. 

Let $A= \otimes^\infty M_k$ with connecting morphisms $a \mapsto a \otimes 1$ be a UHF-algebra of type $k^\infty$ and let $H_{\tr}$ be the GNS-Hilbert space associated to the unique trace $\tr$; $H_{\tr}$ is an inductive limit of finite-dimensional Hilbert spaces $H_n$ associated to $M_{k^n}$, $n \in \mathbb{N}$. Let $(A, H_{\tr} , D)$ be the Christensen-Ivan spectral triple associated to an increasing sequence of real numbers $\{\lambda_n\}$. If $v \in H_{\tr}$ is a unit vector we denote by $\phi_v$ the associated vector state on $A$. Denote $A_n = M_{k^n} \simeq M_k ^{\otimes n}$.

\begin{lem}
\label{lem6.1}
Keep the above notation and let $n \in \mathbb{N}$, $v \in H_n$ with $\tr(v^* v) = 1$. Then
\begin{equation*}
d_D(\phi_v , \tr)= \sup_{x \in A_n, \; \|[D,x]\| \leq 1} | \tr (x) - \phi_v (x)|.
\end{equation*}
\end{lem}
\proof It is enough to show that for every $\tilde{x} \in \bigcup_h A_h$ with $\|[D,\tilde{x}]\| \leq 1$ there is $x \in A_n$ such that $\|[D,x]\| \leq 1$ and $|\tr  (x) - \phi_v (x) | = |\tr(\tilde{x})-\phi_v (\tilde{x})|$. Hence let $\tilde{x} \in \bigcup_h A_h$ and consider the conditional expectation $E_n : \bigcup_h A_h \rightarrow A_n$ given by $E_n (\sum x_1 \otimes x_2 \otimes ... \otimes x_n \otimes x_{n+1} \otimes...)= \sum  x_1 \otimes x_2 \otimes ... \otimes x_n \otimes \tr (x_{n+1} ) \otimes \tr (x_{n+2}) \otimes ...$. We claim that $E_n (\tilde{x}) \in A_n$ is the desired element. For note that for every $v_1, v_2 \in H_n$ we have 
\begin{equation}
\label{eqeq}
\braket{v_1, E_n (\tilde{x}) v_2} = \langle v_1, \tilde{x} v_2\rangle
\end{equation}
 and so, denoting by $P_n$ the orthogonal projection from $H_{\tr}$ onto $H_n$, we have
\begin{equation*}
P_n \tilde{x} P_n = P_n E_n (\tilde{x}) P_n.
\end{equation*}
 Moreover, $E_n (\tilde{x})$ leaves every $H_k$ invariant for $k \geq n$. Hence
 \begin{equation*}
 [D, E_n (\tilde{x})] = \sum_k \alpha_k (Q_k E_n (\tilde{x}) - E_n (\tilde{x}) Q_k) = \sum_{k=1}^n \alpha_k P_n [Q_k, E_n (\tilde{x})] P_n = P_n [D,\tilde{x}] P_n,
 \end{equation*}
 from which we obtain $\| [D,E_n (\tilde{x})]\| \leq \|[D,\tilde{x}]\|$. Furthermore, it follows again from (\ref{eqeq}) that $|\tr (\tilde{x})- \phi_v (\tilde{x})|= |\tr(E_n (\tilde{x}))- \phi_v (E_n (\tilde{x}))|$. \hfill $\Box$
 
 The above Proposition is the key observation for the explicit computation of the distance between the trace and certain states in the case of the CAR-algebra. This computation reveals the impossibility for the automorphisms $\sigma_{i,j}$ (see the discussion before Lemma \ref{leml}) to belong to the $\ISO$-group of a Christensen-Ivan spectral triple. We recall that $M_2$ is linearly generated by the Pauli matrices
 \begin{equation*}
	\sigma_1 = \left(\begin{array}{cc}	0	&	1	\\
								1	&	0	\end{array}\right),\quad
								\sigma_2 = \left(\begin{array}{cc}	0	&	-i	\\
								i	&	0	\end{array}\right),\quad	
 \sigma_3 = \left(\begin{array}{cc}	1	&	0	\\
								0	&	-1	\end{array}\right),\quad
								\sigma_4 = \left(\begin{array}{cc}	1	&	0	\\
								0	&	1	\end{array}\right).
								\end{equation*}
 
In order to clarify the next result, we first illustrate a concrete basic example in which the distance of the trace from a specific vector state in the case $A$ is the CAR-algebra can be computed. In the following we fix a Christensen-Ivan spectral triple associated to a sequence of diverging pairwise different eigenvalues for such algebra. Let
\begin{equation*}
v=\left( \begin{array}{cc}	0	&	\sqrt{2}	\\
	0	&	0	\end{array}\right) \quad \in M_2 \otimes 1 
	\end{equation*} 
	considered as a norm one vector in the $GNS$ representation $H_{\tr}$.  
We have $\phi_v (\sigma_1)=\phi_v (\sigma_2)=0$ and $\phi_v (\sigma_3)=\phi_v (\sigma_4)=1$. We want to show that the $\sup$ appearing in the definition of the distance between $\phi_v$ and $\tr$ can be computed on multiples of $\sigma_3$; in order to do so we need to show that for every $\alpha_1 = \gamma, \alpha_2 = \delta, \alpha_3 = \beta , \alpha_4 = \alpha \in \mathbb{C}$ we have $\|[D, \sum_{i=1}^4 \alpha_i \sigma_i]\| \geq \| [D, \alpha_3 \sigma_3]\|$. Let  
\begin{equation*}
x=  \left(\begin{array}{cc}	\alpha + \beta 	&	\gamma - i \delta	\\
								\gamma + i\delta	&	\alpha - \beta	\end{array}\right) \quad \in M_2.
								\end{equation*}
Using the $C^*$-identity we see that the norm of the commutator $[D,x]$ is given by $\|[D,x]\|= \lambda_1 \max\{\|x_{0,1}\|, \|x_{1,0}\|\}$, where $x_{0,1} =P_0 x (P_1-P_0)$ and $x_{1,0} = (P_1 - P_0) x P_0$. A computation shows that for $w \in \mathbb{C}$ we have
\begin{equation*}
x_{1,0} w = \left( \begin{array}{cc} \beta w 	& 	(\gamma - i \delta) w	\\
					(\gamma + i \delta) w 	&	-\beta w	\end{array}\right),
\end{equation*}
from which we obtain
\begin{equation*}
\|x_{1,0}\|^2 = |\beta|^2 + \frac{1}{2} (|\gamma + i \delta|^2 + |\gamma - i\delta|^2) \geq |\beta|^2  .
\end{equation*}	
Note now that $(\beta\sigma_3)_{0,1} : \mathbb{C}\xi_{\tr}^{\perp} \subset H_1 \rightarrow \mathbb{C}$ is given by
\begin{equation*}
(\beta \sigma_3)_{0,1}  \left(\begin{array}{cc}	w_{1,1}	&	w_{1,2}	\\
								w_{2,1}	&	-w_{1,1}	\end{array}\right) = \beta w_{1,1}
\end{equation*}
and so $\|(\beta  \sigma_3)_{0,1}\| = |\beta|$. Similarly we have $\| (\beta \sigma_3)_{1,0}\|=|\beta|$. In particular $\|x_{1,0}\| \geq \max\{\| (\beta \sigma_3)_{1,0}\|, \| (\beta \sigma_3)_{0,1}\|\}$ and so $\|[D,x]\| \geq \| [D, \beta \sigma_3]\|$.
It follows from Lemma \ref{lem6.1} that
\begin{equation*}
d_D (\phi_v , \tr)= \sup_{\beta \in \mathbb{C} \; : \; \| [D,\beta\sigma_3]\| \leq 1} |\tr (\beta\sigma_3) - \phi_v (\beta \sigma_3)| =\sup_{\beta \in \mathbb{C} \; : \; \| [D,\beta\sigma_3]\| \leq 1} |\beta| .
\end{equation*}
In order to obtain this value we only need the value of $\|[D, \beta \sigma_3]\|$, which we already computed and is given by $|\beta| \lambda_1$ and so 
\begin{equation*}
\label{eqda1}
d_D (\phi_v , \tr) = \lambda_1^{-1} \ . 
\end{equation*}
 
 This specific example already contains all the ideas needed for the general case.
 
  \begin{prop}
Let $A=\otimes^\infty M_2$, for $l \in \{1,2,3\}$ denote $B_l:=   \{ v \in M_2 \; | \; \; \phi_v (\sigma_l) = \phi_v (\sigma_4) =1, \; \phi_v (\sigma_j) =0 \mbox{ for } j \in \{1,2,3\} \backslash \{l\} \}  \subset M_2$ ($\phi_v$ is the vector state associated to $v$: $\phi_v (\cdot)= \tr (v^* (\cdot) v)$) and let $v \in \ad(U(M_2))(B_l)$ for some $l$. For every $n \in \nn$ we have
\begin{equation*}
d(\phi_{1^{\otimes n} \otimes v}, \tr) = \frac{1}{\lambda_{n+1}}.
\end{equation*}
The supremum for the distance is attained on the element $1^{\otimes n}\otimes \sigma_l /\| [D,1^{\otimes n} \otimes \sigma_l]\|  \in 1^{\otimes n} \otimes M_2 \subset A_{n+1}$.  
\end{prop}
\proof 
First of all note that we can suppose $v \in B_l$, since $\ad(1^{\otimes n} \otimes  U(M_2)) \subset \Iso (A, D, H_{\tr})$ in virtue of Theorem \ref{p1}. Let $\sigma_i$, $i=1,...,4$ be the Pauli matrices and $v$ be as in the statement. For every 
\begin{equation*}
x= \sum_{i_1,i_2,...,i_{n+1} =1}^4 \alpha_{i_1,i_2,...,i_{n+1}} \sigma_{i_1} \otimes \sigma_{i_2} \otimes ... \otimes \sigma_{i_{n+1}}\in A_{n+1}
\end{equation*} we denote 
\begin{equation*}
\tilde{x}=x - \sum_{i_1,i_2,...,i_n =1}^4 \alpha_{i_1,i_2,..,i_n, 4} \sigma_{i_1} \otimes \sigma_{i_2} \otimes ... \otimes \sigma_{i_n} \otimes 1 \qquad \in A_{n+1}. \end{equation*}
 We want to show that
\begin{equation*}
\begin{split}
d(\phi_{1^{\otimes n} \otimes v}, \tr) &= \sup_{x \in A_{n+1}, \| [D,x]\| \leq 1} |\phi_{1^{\otimes n} \otimes v} (x) - \tr (x)| = \sup_{x \in A_{n+1}, \|[D,\tilde{x}]\|\leq 1} |\phi_{1^{\otimes n} \otimes v }(\tilde{x})|\\
&= \sup_{x \in A_{n+1}, \|[D,\tilde{x}]\|\leq 1} |\alpha_{4,4,...,4,l}|. \end{split}
\end{equation*}
 For note that $\phi_{1^{\otimes n} \otimes v} (\sigma_{i_1} \otimes \sigma_{i_2} \otimes ... \otimes \sigma_{i_{n+1}}) = \prod_{k=1}^n\delta_{i_k,4} (\delta_{i_{n+1},l} + \delta_{i_{n+1},4})$ and $\tr(\sigma_{i_1} \otimes \sigma_{i_2} \otimes ... \otimes \sigma_{i_{n+1}}) = \prod_{k=1}^{n+1}\delta_{i_k, 4}$. Hence for every $x=\sum_{i_1,i_2,...,i_{n+1} =1}^4 \alpha_{i_1,i_2,...,i_{n+1}} \sigma_{i_1} \otimes \sigma_{i_2} \otimes ... \otimes \sigma_{i_{n+1}}$ we have $|\phi_{1^{\otimes n} \otimes v} (x) - \tr (x)| = |\alpha_{4,4,...,4,l}| $. We are then left to show that for every $x \in A_{n+1}$, $\|[D,x]\| \geq \|[D,\tilde{x}]\|$. In order to do so, let $(i_1,i_2,...,i_{n+1})\in \{1,2,3,4\}^{\times n} \times \{1,2,3\}$ and $w \in H_n$. Then $P_n \sigma_{i_1} \otimes \sigma_{i_2} \otimes ... \otimes \sigma_{i_{n+1}} w= \sigma_{i_1} \otimes \sigma_{i_2} \otimes ... \otimes \sigma_{i_n} w \otimes \tr (\sigma_{i_{n+1}}) =0$ and so $P_n \sigma_{i_1} \otimes \sigma_{i_2} \otimes ... \otimes \sigma_{i_{n+1}}P_n=0$, which gives $P_n [D, \sigma_{i_1} \otimes \sigma_{i_2} \otimes ... \otimes \sigma_{i_{n+1}}] P_n=0$ so that $P_n [D,\tilde{x}]P_n=0$. On the other hand, $x-\tilde{x} \in A_n$ commutes with $Q_k$ for every $k >n$; it follows that $ [D,x-\tilde{x}]Q_k=Q_k [D,x-\tilde{x}] =0$ for each such $k$. Also, a direct computation shows
 \begin{equation*}
 \begin{split}
 Q_{n+1}[D, \tilde{x}]Q_{n+1} &= Q_{n+1} (\sum_i \alpha_i Q_i \tilde{x} - \sum_i \alpha_i \tilde{x} Q_i )Q_{n+1}\\ & = \alpha_{n+1} Q_{n+1} \tilde{x} Q_{n+1} - \alpha_{n+1} Q_{n+1} \tilde{x} Q_{n+1} =0.
 \end{split}
 \end{equation*}
 Thus $[D,\tilde{x}]$ and $[D, x-\tilde{x}]$ satisfy the following properties:
 \begin{equation*}
 P_n [D, \tilde{x}]P_n =0, \qquad P_{n+1} [D, \tilde{x}] P_{n+1} =[D, \tilde{x}], \qquad Q_{n+1} [D, \tilde{x}] Q_{n+1} =0,
 \end{equation*}
 \begin{equation*}
 [D,x-\tilde{x}] Q_k = Q_k [D, x-\tilde{x}] =0 \qquad \forall k>n.
 \end{equation*}
 
 Hence we want to check that for every $y,z \in \mathbb{B}(H_{\tr})$ with $P_{n+1} y P_{n+1} = y$, $P_n y P_n =Q_{n+1} y Q_{n+1}=0$ and $ zQ_k=Q_k z =0$ for every $k >n$ we have $\|y\| \leq \|y+z\|$. Under these assumptions we have $\| y|_{K_i}\| \leq \| (y+z)|_{K_i}\|$ for every $i=0,1,...,n+1$. Now, $ y = P_n y Q_{n+1} + Q_{n+1} y P_n $ and so 
 \begin{equation*} 
 \begin{split}
 y^* y&= (Q_{n+1} y^* P_n + P_n y^* Q_{n+1})(P_n y Q_{n+1} + Q_{n+1} y P_n)\\&=Q_{n+1} y^* P_n y Q_{n+1} + P_n y^* Q_{n+1} y P_n\end{split}
 \end{equation*} so that $\|y\| = \max \{ \|P_n y Q_{n+1}\|, \|Q_{n+1} y P_n\|\}$. We have $P_n y Q_{n+1} = y|_{K_{n+1}}$ and so $\| P_n  y Q_{n+1}\| \leq \|y+z\|$. Now note that the above conditions for $y$ and $z$ are preserved under taking adjoints and so $\|y^*|_{K_{n+1} }\| \leq \|(y+z)^*|_{K_{n+1} }\| \leq \| (y+z)^*\| = \|y+z\|$. But again $y^*|_{K_{n+1}} = P_n y^* Q_{n+1} = (Q_{n+1} y P_n)^*$, which entails $\|y^*|_{K_{n+1}}\| = \|Q_{n+1} y P_n\|$.\\
 The next step is to show that for every $x \in A_{n+1}$ as above we have $|\alpha_{4,4,...,4,l}| \leq \| P_0 \tilde{x} Q_{n+1}\| \leq  \lambda_{n+1} ^{-1}$. The first inequality is obtained by evaluating $P_0 \tilde{x} Q_{n+1}$ on the unit vector $1^{\otimes n} \otimes \sigma_l \in K_{n+1}$, which gives, using the orthogonality relations for the Pauli matrices, $P_0 \tilde{x} (1^{\otimes n} \otimes \sigma_l) =  \alpha_{4,4,...,4,l}$. For what concerns the second inequality, note that 
 \begin{equation*}
 \begin{split}
 [D,\tilde{x}]&= \lambda_{n+1} Q_{n+1} \tilde{x} P_0 + (\lambda_{n+1} - \lambda_1) Q_{n+1} \tilde{x} Q_1 +...  + (\lambda_{n+1} - \lambda_n) Q_{n+1} \tilde{x} Q_n -\\ & \lambda_{n+1} P_0 \tilde{x} Q_{n+1} - (\lambda_{n+1} - \lambda_1) Q_1 \tilde{x} Q_{n+1} -...- (\lambda_{n+1} - \lambda_n ) Q_n \tilde{x} Q_{n+1}.\end{split}\end{equation*} Since $\|[D,\tilde{x}]\| \leq 1$ we get $\| P_0 \tilde{x} Q_{n+1}\| \leq \lambda_{n+1}^{-1}$.\\
  We are left to prove that 
  \begin{equation*} |\phi_{1^{\otimes n}\otimes v} (1^{\otimes n} \otimes \sigma_l /\|[D, 1^{\otimes n}\otimes \sigma_l]\|) - \tr(1^{\otimes n} \otimes \sigma_l /\|[D, 1^{\otimes n}\otimes \sigma_l]\|)| = \lambda_{n+1}^{-1} . 
  \end{equation*}This follows from the fact that $\|[D,1^{\otimes n}\otimes \sigma_l]\| = \lambda_{n+1}$. For note that 
  \begin{equation*}
  \begin{split}
  \|[D, 1^{\otimes n}\otimes \sigma_l]\|  & = \max\{\|\lambda_{n+1} P_0(1^{\otimes n}\otimes \sigma_l )Q_{n+1} + (\lambda_{n+1} - \lambda_n) Q_1 (1^{\otimes n}\otimes \sigma_l )Q_{n+1} +...\\ &\quad + (\lambda_{n+1} - \lambda_n) Q_n (1^{\otimes n}\otimes \sigma_l) Q_{n+1}\|, \|\lambda_{n+1} Q_{n+1} (1^{\otimes n}\otimes \sigma_l) P_0 \\ &\quad + (\lambda_{n+1} - \lambda_1) Q_{n+1}(1^{\otimes n}\otimes \sigma_l) Q_1 +... + (\lambda_{n+1} - \lambda_n) Q_{n+1} (1^{\otimes n}\otimes \sigma_l) Q_n\|\} \\ &=  \|\lambda_{n+1} P_0(1^{\otimes n}\otimes \sigma_l )Q_{n+1} + (\lambda_{n+1} - \lambda_n) Q_1 (1^{\otimes n}\otimes \sigma_l )Q_{n+1} +...\\ & \quad + (\lambda_{n+1} - \lambda_n) Q_n (1^{\otimes n}\otimes \sigma_l) Q_{n+1}\|\\ &=\|\lambda_{n+1} Q_{n+1} (1^{\otimes n}\otimes \sigma_l) P_0 + (\lambda_{n+1} - \lambda_1) Q_{n+1}(1^{\otimes n}\otimes \sigma_l) Q_1 +... \\ & \quad + (\lambda_{n+1} - \lambda_n) Q_{n+1} (1^{\otimes n}\otimes \sigma_l) Q_n\| .\end{split}\end{equation*} Let $\alpha \oplus v_1 \oplus v_2 \oplus ... \oplus v_n \in H_0 \oplus K_1 \oplus ... \oplus K_n$ with $|\alpha|^2 +\sum_{i=1}^n \| v_i\|^2=1$. We have 
  \begin{equation*}
  \begin{split}
 & \|  ( \lambda_{n+1} Q_{n+1} (1^{\otimes n}\otimes \sigma_l) P_0 + (\lambda_{n+1} - \lambda_1) Q_{n+1}(1^{\otimes n}\otimes \sigma_l) Q_1 +... \\ \qquad & + (\lambda_{n+1} - \lambda_n) Q_{n+1} (1^{\otimes n}\otimes \sigma_l) Q_n)(\alpha \oplus v_1 \oplus ... \oplus v_n)\|^2\\ & =\lambda_{n+1}^2 |\alpha|^2 + (\lambda_{n+1}-\lambda_1)^2 \|v_1\|^2 +... + (\lambda_{n+1} - \lambda_n)^2 \|v_n\|^2\\ & \leq \lambda_{n+1}^2 (|\alpha|^2+\sum_{i=1}^n \| v_i\|^2) = \lambda_{n+1}^2 , \end{split}\end{equation*}which gives $\| \lambda_{n+1} Q_{n+1} (1^{\otimes n}\otimes \sigma_l) P_0 + (\lambda_{n+1} - \lambda_1) Q_{n+1}(1^{\otimes n}\otimes \sigma_l ) Q_1 +... + (\lambda_{n+1} - \lambda_n) Q_{n+1} (1^{\otimes n}\otimes \sigma_l) Q_n\| =\lambda_{n+1}$. \hfill $\Box$
 
 \begin{cor}
 Let $A= \otimes^\infty M_2$ and suppose that the eigenvalues $\lambda_n$ of the Dirac operator are pairwise distinct. Let $\alpha \in \Aut (A)$ be given with the following property: there exist $k\neq   m  \in \nn$, $l,i \in \{1,2,3\}$, $v \in \ad(U(M_2))(B_l)$ such that $\alpha (1^{\otimes m} \otimes v) \in 1^{\otimes k} \otimes \ad(U(M_2))(B_i)$. Then $\alpha $ does not belong to $\ISO (A,H_{\tr}, D)$. In particular this applies to the automorphisms $\sigma_{i,j}$ given by permutation of the $i$-th and $j$-th tensor factors.
 
 \end{cor}
 
 In the end of this section we include a related result the authors think is of independent interest.
 
\begin{prop}
 Denote by $\varphi : \bigcup_n A_n \rightarrow \bigcup_n A_n$ the endomorphism given by $\varphi (y) = 1 \otimes y$. For every $n \in \nn$, $c(n) >0$ such that $(c(n)+1)\lambda_{n } \leq \lambda_{n+1} $ we have  
  \begin{equation*}
  \|[D, x]\| \leq c(n)^{-1} \lambda_{n }^{-1} \lambda_1\| [D, \varphi^n (x)]\| 
  \end{equation*}
  for every $x \in A_1$. In particular if $2\lambda_1< \lambda_2$, then $\|[D,\varphi(x)]\| \geq \frac{\lambda_2-\lambda_1}{\lambda_1} \|[D, x ]\|$ for very $x \in A_1$, with $\frac{\lambda_2-\lambda_1}{\lambda_1}>1$.
  \end{prop}
  \proof Let $x \in A_1$ and $\xi$ be the cyclic vector. First we show that $\|[D,x]\| \leq \lambda_1 \|Q_1x \xi\|_2$. For note that the $C^*$-identity gives $\|[D,x]\| =\lambda_1 \max\{ \| P_0 x Q_1\|, \|Q_1 xP_0\|\}$. Now, for every $v \in Q_1 H$ we have $P_0 xv= \tr (xa) \xi= \tr ((x-\tr (x))a) \xi$, where $a\xi =v$; hence $\|P_0 x Q_1\| \leq\|Q_1x\xi\|_2$. Similarly, for $\alpha \in \mathbb{C}$, $Q_1 x \alpha\xi  = \alpha (x-\tr(x))\xi$ and so again $\|Q_1 x P_0\| \leq \|Q_1x \xi\|_2$.\\
 A direct computation gives $P_n[D, \varphi^n (x)]P_n=0$; hence $$P_n [D, \varphi^n (x)]^* [D,\varphi^n (x)] P_n= P_n[D,\varphi^n (x)]^* Q_{n+1} [D,\varphi^n (x)] P_n$$ and we only need to show that $\| P_n[D,\varphi^n (x)]^* Q_{n+1} [D,\varphi^n (x)] P_n\| \geq c(n)^2 \lambda_n^2 \|Q_1 x\xi\|^2_2$. The operator $P_n[D,\varphi^n (x)]^* Q_{n+1} [D,\varphi^n (x)] P_n$ is explicitly given by 
\begin{equation*}
\begin{split}
 P_n & [ D,\varphi^n (x)]^*  Q_{n+1} [D,\varphi^n (x)] P_n  = \sum_{i,j=0}^n (\lambda_{n+1} - \lambda_i)(\lambda_{n+1} - \lambda_j) Q_i (1^{\otimes n} \otimes x^* )Q_{n+1} (1^{\otimes n} \otimes x )Q_j\\
 &= \sum_{i=0}^n (\lambda_{n+1} - \lambda_i)^2 Q_i( 1^{\otimes n}\otimes x^*) Q_{n+1}( 1^{\otimes n} \otimes x )Q_i \\ &+ \sum_{i<j} (\lambda_{n+1} - \lambda_i) (\lambda_{n+1} - \lambda_j) [Q_i (1^{\otimes n} \otimes x^* )Q_{n+1} (1^{\otimes n} \otimes x )Q_j + Q_j (1^{\otimes n} \otimes x^* )Q_{n+1} (1^{\otimes n} \otimes x) Q_i]\\
& \geq c(n)^2 \lambda_n^2 (P_n (1^{\otimes n} \otimes x^* )Q_{n+1} (1^{\otimes n} \otimes x )P_n).
\end{split}
\end{equation*}
Explicit computations give, for $v=a\xi \in H_n$, $P_n (1^{\otimes n} \otimes x^*) Q_{n+1}  (1^{\otimes n} \otimes x	)P_n (v) = a \otimes\tr(x^* (x-\tr (x))) \xi = a \otimes \| Q_1 x\xi\|_2^2$, from which we obtain $\|P_n (1^{\otimes n}\otimes x)^* Q_{n+1} (1^{\otimes n}\otimes x) P_n\| \geq c(n)^{ 2} \lambda_n^{ 2}\| Q_1 x\xi\|^2_2$ and
\begin{equation*}
\|[D,\varphi^n (x)]\|^2 \geq \| P_n [D, \varphi^n (x)]^* [D,\varphi^n (x)]P_n\| \geq c(n)^2 \lambda_n^2 \|Q_1 x\xi \|_2^2 \geq c(n)^2 \lambda_n^{2} \lambda_1^{-2}\|[D,x]\|^2. \qquad \Box
\end{equation*}

\subsection{The case of the Cantor set}
Let $X$ be a (compact) Cantor set. There is a nested partition $\mathcal{P}=\{\mathcal{P}_i\}$ with $\mathcal{P}_i = \{U_{i,j}, \; j \in  \mathbb{Z}_2^i\}$ of $X$ consisting of clopen sets such that for every $i \in \nn$, $j \in \mathbb{Z}_2^i$ we have $U_{i,j} = U_{i+1, j\oplus 0} \sqcup U_{i+1, j \oplus 1}$ and every continuous function on $X$ is the uniform limit of $\mathcal{P}$-simple functions. If we identify $X$ with $\prod_{j=1}^\infty \{0,1\}$ (with the product topology), an explicit description of such a nested partition can be given in the following way: for every $l \in \mathbb{N}$ let $\pi_l $ be the projection on the factor $X_l= \prod_{j=0}^l \{0,1\}$ of $\prod_{j=0}^\infty \{0,1\}$; given a finite word $w \in X_l$ denote by $U_w$ the clopen set given by $U_w := \{ x \in X \; | \; \pi_l (x) =w\}$ and define $\mathcal{P}_i = \{ U_w \; | \; w \in X_i\}$. The $C^*$-algebra of continuous functions on $X$ is the inductive limit of the sequence of finite dimensional $C^*$-algebras $\cc^{2^i}$ consisting of $\mathcal{P}_i$-simple functions, where the connecting morphism $\cc^{2^i} \rightarrow \cc^{2^{i+1}}$ is given by considering a $\mathcal{P}_i$-simple function as a $\mathcal{P}_{i+1}$-simple function. This gives a concrete realization of $C(X)$ as the $AF$-algebra $\lim_i \cc^{2^i}$ with connecting morphisms $\phi_{i,i+1} (\oplus_{k=1}^{2^i} a_k)= \oplus_{k=1}^{2^i} (a_k \oplus a_k)$. Fix such a nested partition $\mathcal{P}$.
\begin{prop}
\label{prop3.1}
Let $X$ be a (compact) Cantor set and let $\mu$ be the uniform measure on $X$, viewed as a faithful state on $C(X)$. Let $(C(X), H_\mu,D)$ be the Christenen-Ivan spectral triple associated with a non-degenerate sequence $\{\lambda_n\}$ and to the $AF$-structure induced by $\mathcal{P}$. Then, denoting by $S_2 \simeq \mathbb{Z}_2$ the group of permutations of the set $\{0,1\}$, we have
\begin{equation*}
\Iso (C(X),  H_\mu , D) = \ltimes_{i=1}^\infty S_2^{2^{i-1}}
\end{equation*}
\end{prop}
\proof
Let $\alpha \in \Iso(C(X),H_\mu,D)$. Since $\alpha$ respects the filtration of the algebra (by Theorem \ref{p1}), for every $i \in \nn$ it induces an automorphism on $\cc^{2^i}$ which leaves $\oplus^{2^{i-1}} (\diag (\cc))$ invariant; by Gelfand duality $\alpha|_{\cc^{2^i}}$ is an element of $S_2^{2^{i-1}} \times G_i$, where $G_i$ is the group of automorphisms of $\oplus^{2^{i-1}} (\diag (\cc))$ which respects the filtration up to the $(i-1)$-th step. Hence $\alpha|_{\cc^{2^i}}$ belongs to $\prod_{k=1}^i S_2^{2^{k-1}}$. Clearly every element of $\prod_{i=1}^\infty S_2^{2^{i-1}}$ induces an automorphism $\alpha$ of $C(X)$ which belongs to $\Iso(C(X),H_\mu,D)$. The group structure follows at once. \hfill $\Box$

\begin{prop}
\label{prop3.2}
Let $X$ be a (compact) Cantor set and consider the same setting as in Proposition \ref{prop3.1}, in the particular case of a Dirac operator associated to a sequence of real numbers $\{\lambda_n\}$ with $\lambda_n=\gamma^{-n+1}$, where $0<\gamma< (3-\sqrt{5})/2$. Then
\begin{equation*}
\ISO(C(X),  H_\mu, D) = \Iso (C(X),  H_\mu , D).
\end{equation*}
\end{prop}
\proof
Let $\alpha$ be an automorphism of $C(X)$ (by an abuse of notation we will eventually consider $\alpha$ as a homeomorphism of $X$) which leaves the metric $d_\gamma$ induced by $D$ invariant. Following \cite{AnCh2}, for $x, y \in X$ let $m(x,y)$ be the least integer $n$ satisfying $x(n) \neq y(n)$, where we view $x$ and $y$ as elements of $\prod \mathbb{Z}_2$. By Theorem 4.1 (ii) of \cite{AnCh2}, under the choice of $\gamma$ as in the statement, if $x,y \in X$ satisfy $m(\alpha x, \alpha y) \neq m(x,y)$, then $d_\gamma (\alpha x, \alpha y) \neq d_\gamma (x,y)$; hence, if $\alpha$ leaves $d_\gamma$ invariant, then it also leaves $m$ invariant. We will see that this condition guarantees $\alpha(A_n) \subset A_n$ for every $n$. First of all note that, since the image under $\alpha$ of a clopen set is a clopen set, for every $n \in \nn$ the image of a minimal projection $p \in A_n$ under $\alpha$ is a finite sum of minimal projections in $\bigcup_i A_i$. Suppose now that there are integers $n<k$ such that the image of $\alpha|_{A_n}$ contains an element $b$ of $A_k$ which is not in $A_n$. If we write $b=\oplus_{i=1}^{2^n}(\oplus_{j=1}^{2^{k-n}} a_{i,j})$, the condition that $b$ does not belong to $A_n$ is equivalent to the existence of $i \in \{1,...,2^n\}$ and $j_1, j_2 \in \{1,...,2^{k-n}\}$ such that $a_{i,j_1} \neq a_{i,j_2}$; it follows that there are two characters $x$ and $y$ with $m(x,y) >n$ and such that $x|_{A_n} = y|_{A_n}$ and $x \circ \alpha |_{A_n} \neq y \circ \alpha |_{A_n}$, which gives $m (\alpha x, \alpha y) \leq n$. \hfill $\Box$

\section{Isometries of crossed products}
In this section we show that, under suitable hypothesis, given an action of a discrete group on a $C^*$-algebra, the crossed product automorphisms given in \cite{DaPa} are elements of the $\Iso$-group of the crossed product spectral triple introduced in \cite{HSWZ}. 

\smallskip 
Let $A$ be a unital $C^*$-algebra, $G$ a countable discrete group and $\alpha : G \ra \Aut (A)$ an action. As in \cite{DaPa} let $Z^1 (G, U(A)):= \{ c : G \rightarrow U(A) \; |\; c(g h)= c(g) \alpha_g (c(h)) \; \forall g,h \in G\}$. We consider the set $\mathcal{G}(A,\alpha) := \{ (c, \beta, \sigma) \; | \; c \in Z^1 (G , U(A)), \beta \in \Aut(A), \sigma \in \Aut(G), \; \beta \circ \alpha_g = \ad c_{\sigma (g)} \circ \alpha_{\sigma (g)} \circ \beta \} \subset Z^1 (G, U(A)) \times \Aut(A) \times \Aut(G)$, endowed with the group structure $(c,\beta , \sigma) (c', \beta' , \sigma')=( \beta (c' \circ \sigma^{-1}) \cdot c , \beta \circ \beta', \sigma \circ \sigma')$, $(c,\beta,\sigma)^{-1} = (\beta^{-1} ( c^*\circ \sigma), \beta^{-1}, \sigma^{-1})$ (see \cite{DaPa}). There is a group homomorphism $\Phi : \mathcal{G}(A,\alpha) \rightarrow \Aut(A \rtimes_{\alpha, r} G, A)$ explicitely given by $\Phi (c,\beta,\sigma)(\sum_{g \in G} a_g \lambda_g) = \sum_{g \in G} \beta (a_g) c_{\sigma(g)} \lambda_{\sigma (g)}$ (the fact that this is actually an automorphism of the reduced crossed product follows since it is continuous in the $L^1$-norm). As shown in \cite{DaPa}, under suitable assumptions $\Phi (\mathcal{G}(A,\alpha)) = \Aut (A\rtimes_r G, A)$.

If $(A,H, D)$ is a spectral triple and $\mathcal{A}$ is a dense $*$-subalgebra of $A$ such that $[D,a]$ extends to a bounded operator on $H$ for every $a \in \mathcal{A}$, we will say that $(\mathcal{A},  H, D)$ is a spectral triple for $A$. This notion coincides with the notion of odd spectral triple considered in \cite{HSWZ} Definition 2.1. In \cite{HSWZ} Theorem 2.7 it is proved that if $(\mathcal{A},H, D)$ is a spectral triple for $A$ and $G$ is a discrete countable group acting on $A$ equipped with a proper translation bounded function $l : G \rightarrow \mathbb{Z}$ (\cite[Example 2.4]{HSWZ}) then, under some additional hypotheses (see the statement of Theorem \ref{thm1} below), $(C_c (G,\mathcal{A}),  (H\otimes l^2(G))^{\oplus 2}, D_l)$ is a spectral triple for $A \rtimes_{r} G$, where $D_l$ is given by
\begin{equation*}
D_l =\left(\begin{array}{cc}	0	&	D\otimes 1 -i \otimes M_l\\
	D \otimes 1 + i \otimes M_l 	&	0
	\end{array}\right)
\end{equation*}
and $M_l$ is the usual self-adjoint extension of the unbounded operator of multiplication by $l$ on $l^2(G)$.

\begin{thm}
\label{thm1}
Let $A$ be a unital $C^*$-algebra endowed with an action $\alpha$ of a countable discrete group $G$, equipped with a proper translation bounded function $l : G \rightarrow \mathbb{Z}$. Suppose we are given a spectral triple $(\mathcal{A},  H, D)$ for $A$ which satisfies the hypothesis of Theorem 2.7 of \cite{HSWZ}, namely: $\alpha_g (\mathcal{A}) \subset \mathcal{A}$ for every $g \in G$, $\sup_{g \in G} \| [D,  \alpha_g (a) ]\| < \infty$ for every $a \in \mathcal{A}$. Let $(c,\beta, \sigma) \in \mathcal{G}(A,\alpha)$ be such that $\beta \in \Iso (A,  H, D)$, $l \circ \sigma = l$ and $c_g \in \mathbb{C} \cdot 1_A$ for every $g \in G$. Then $\Phi (c,\beta,\sigma) \in \Iso (A\rtimes_{\alpha,r} G,  (H \otimes \ell^2 (G))^{\oplus2}, D_l)$, where $D_l$ is the above Dirac operator.
\end{thm}
\proof 
Let $U_\beta \in U(H)$ be a unitary implementing the automorphism $\beta$. Let $U_{(c\beta,\sigma)}$ be the unitary operator in $U(H \otimes \ell^2 (G))$ given, for $\xi \in H $ and $g \in G$, by $U_{(c,\beta,\sigma)}(\xi \otimes \delta_g)=(c_{\sigma (g^{-1})}^* U_\beta \xi)\otimes \delta_{\sigma (g)}$. For $a \in A$, $\xi \in H$ and $g, h \in G$ we compute
\begin{equation*}
\begin{split}
\Phi(c,\beta,\sigma) (a\lambda_g)& U_{(c,\beta,\sigma)} \xi \otimes \delta_h = (\beta(a) c_{\sigma (g)} \lambda_{\sigma(g)} )((c_{\sigma(h^{-1})}^* U_\beta \xi)\otimes \delta_{\sigma (h)})\\
&=(\alpha_{\sigma (gh)^{-1}}(\beta(a) c_{\sigma (g)}) c_{\sigma(h^{-1})}^* U_\beta \xi)\otimes \delta_{\sigma (gh)}\\
&=(c_{\sigma(gh)^{-1}}^* \beta (\alpha_{(gh)^{-1}} (a)) c_{\sigma (gh)^{-1}}\alpha_{\sigma(gh)^{-1}} (c_{\sigma (g)})  c_{\sigma (h)^{-1}}^* U_\beta \xi)\otimes \delta_{\sigma (gh)}\\
&=(c_{\sigma(gh)^{-1}}^* \beta (\alpha_{(gh)^{-1}} (a)) c_{\sigma(h)^{-1}} c_{\sigma(h)^{-1}}^* U_\beta \xi)\otimes \delta_{\sigma (gh)}\\
&=(c_{\sigma(gh)^{-1}}^* \beta (\alpha_{(gh)^{-1}} (a)) U_\beta \xi)\otimes \delta_{\sigma (gh)}.
\end{split}
\end{equation*}
On the other hand
\begin{equation*}
\begin{split}
U_{(c,\beta,\sigma)}& a\lambda_g (\xi \otimes \delta_h)=U_{(c,\beta,\sigma)} (\alpha_{(gh)^{-1}} (a) \xi) \otimes \delta_{gh}\\
&= (c_{\sigma(gh)^{-1}}^* U_\beta \alpha_{(gh)^{-1}} (a) \xi) \otimes \delta_{\sigma (gh)}= (c_{\sigma(gh)^{-1}}^* \beta (\alpha_{(gh)^{-1}} (a) )U_\beta\xi) \otimes \delta_{\sigma (gh)}.
\end{split}
\end{equation*}
Hence $U_{(c,\beta,\sigma)}$ implements the automorphism $\Phi(c,\beta,\sigma)$ in $H \otimes \ell^2 (G)$. Note that in the above computation we did not use the fact that $c$ is a character. Represent now $A \rtimes_{\alpha, r} G$ diagonally on $(H \otimes \ell^2 (G))^{\oplus2})$; then $\Phi(c,\beta,\sigma)$ is unitarily implemented by $U_{(c,\beta,\sigma)} \otimes 1_2$ in this representation. Let $\{\lambda_k\}_{k \in \nn}$ be the eigenvalues of $D$ and $\{P_k\}$ the corresponding spectral projections. Under the identification $H \otimes \ell^2 (G) \simeq \ell^2 (G, H)$, we can write $\xi \in (H \otimes \ell^2 (G))^{\oplus2}$ as $(\xi_g)_{g \in G}$, with $\xi_g = (\xi_1 , \xi_2)_g \in H \oplus H$ for every $g \in G$ and the domain of the unbounded operator $D_l$ is given by the set of vectors $\xi=(\xi_1 , \xi_2) \in (H \otimes l^{2} (G))^{\oplus2}$ which satisfy
\begin{equation*}
 \sum_{g \in G} \sum_{k \in \nn} (\lambda_k^2 + l(g)^2) \frac{1}{4}\left\| \left(\begin{array}{c} P_k \xi_1 + \frac{\lambda_k - i l(g)}{\sqrt{\lambda_k^2 + l(g)^2}} P_k \xi_2 \\
 \frac{\lambda_k + il(g)}{\sqrt{\lambda_k^2 + l(g)^2}} P_k \xi_1 + P_k \xi_2 \end{array}\right)_g \right\|^2 < \infty.
 \end{equation*}
Under the hypothesis that $U_\beta$ commutes with the spectral projections of $D$, that $c$ is a character and that $l\circ \sigma = l$, this relation is left invariant under $U_{(c,\beta,\sigma)} \otimes 1_2$. In order to have that $[D_l, U_{(c,\beta,\sigma)}]=0$ it is enough to check that
$[D, c_g^* U_\beta]=0$ for every $g \in G$ and that $[1 \otimes M_l, U_{(c,\beta,\sigma)}]=0$; since $c$ is a character, these follow from $[D,U_\beta]=0$ and $l\circ \sigma = l$ respectively. \hfill $\Box$

\begin{oss} 
Let $A$ be a unital $C^*$-algebra endowed with an action of a countable discrete group $G$, $l:G \rightarrow \mathbb{Z}$ a proper translation bounded function and denote $\mathcal{G}_{(l,D,H)} (A, \alpha) := \{ (c, \beta, \sigma) \in \mathcal{G} (A, \alpha) \; | \; l \circ \sigma = l, \; c \in \chara{G},\; \beta \in \Iso (A,H, D)\}$. Let $(\mathcal{A}, H, D)$ be a spectral triple for $A$ in the sense of \cite{HSWZ} Definition 2.1 and suppose that $\alpha_g (\mathcal{A}) \subset \mathcal{A}$ for every $g \in G$, $\sup_{g \in G} \| [D, \pi_\psi (\alpha_g (a))]\| < \infty$ for every $a \in \mathcal{A}$. Theorem \ref{thm1} states that the map $\Phi : \mathcal{G}_{(l,D,H)} (A, \alpha) \ra \Aut (A\rtimes_{\alpha, r} G, A)$ sets up an injective group homomorphism  $\Phi : \mathcal{G}_{(l,D,H)} (A, \alpha) \ra \Iso (A\rtimes_{\alpha, r} G, (H\otimes \ell^2 (G))^{\oplus 2}, D_l)$.
\end{oss} 
 \medskip

\noindent\textbf{Examples:} 
\begin{itemize}
\item[(i)] Let $A$ be a unital $C^*$-algebra endowed with the trivial action of a countable discrete group $G$ admitting a translation bounded function $l: G \rightarrow \mathbb{Z}$ and $(A,H,D)$ be a spectral triple. In this case the crossed product is just the spatial tensor product $A \otimes C^*_r (G)$ and the Dirac operator $D_l$ is the (even) tensor product triple (cfr. the paragraph before Proposition 2.8 in \cite{HSWZ}). Denoting by $\Aut_l (G)$ the subgroup of $\Aut(G)$ given by the automorphisms of $G$ which preserve $l$, the map $\Phi$ gives a bijection between $\chara{G} \times  \Iso(A,H,D) \times \Aut_l (G)$ and $\Iso (A,H,D) \times \Iso (C^*_r (G), M_l, \ell^2 (G))$ (see \cite{LoWu}). This suggests that in general the product of $\Iso$-groups should embed inside the $\Iso$ of some tensor product of spectral triples.
\item[(ii)] Let $A$ be a separable unital $C^*$-algebra and $(\mathcal{A}, H, D)$ a spectral triple for $A$ in the sense of \cite{HSWZ} Definition 2.1. Let also $\alpha \in \Iso (A,H,D)$ be such that $\alpha (\mathcal{A}) \subset \mathcal{A}$ (this is always the case if $A$ is $AF$ and the spectral triple is the one considered by Christensen and Ivan) and $\chi \in S^1$ be a character of $\mathbb{Z}$. Then $\Phi_{(\chi, \alpha, \id)} \in \Iso (A\rtimes_{\alpha, r} \mathbb{Z}, (H\otimes \ell^2(Z))^{\oplus 2}, D_l)$, where $l : \mathbb{Z} \rightarrow\mathbb{Z}$, $l(n)=n$. An example of this situation is given by the odometer actions on the (compact) Cantor set, when we consider the spectral triple of Christensen and Ivan associated to a faithful state. In this case the associated crossed products are the Bunce-Deddens algebras (c.f. \cite{HSWZ}, section 3.1). 
\item[(iii)] More generally, let $(\mathcal{A},H,D)$ be a spectral triple for a unital $C^*$-algebra $A$ in the sense of \cite{HSWZ} Definition 2.1 and consider the subgroup $\Iso (\mathcal{A}, H, D)$ of $\Iso(A,H,D)$ which consists of elements $\alpha$ such that $\alpha (\mathcal{A}) \subset \mathcal{A}$. Then for every countable discrete subgroup $G$ of $\Iso(\mathcal{A},H,D)$ admitting a translation bounded function $l$, every $\alpha$ in the normalizer subgroup of $G$ in $\Iso(A, H,D)$ such that $l=l \circ \ad \alpha$ and character $\chi$ of $G$ we have $\Phi(\chi, \alpha, \ad \alpha) \in \Iso (A \rtimes_{r} G, (H \otimes \ell^2(G))^{\oplus 2}, D_l)$. If $\alpha \in G$, then $\Phi(\chi, \alpha, \ad \alpha)= \chi (\cdot )\ad \lambda_\alpha$, but in general $\lambda_\alpha$ will differ from the unitary $U_{(1, \alpha, \ad \alpha)}$ commuting with the Dirac operator $D_l$. This procedure can be iterated in different ways. 

\item[(iv)] Let $d \in \mathbb{N}$ and $G=\mathbb{Z}^d$ or $\mathbb{F}_d$. Consider the length function on $G$ given by a choice of generators, this is a translation bounded function. Let $\sigma$ be an automorphism of $G$ induced by a permutation of the $d$ generators. Let now $(A,H,D)$ and $\alpha$ be as in example (ii). Considering the action $\tilde{\alpha}$ of $G$ on $A$ induced from the "forgetful" homomorphism $G \rightarrow \mathbb{Z}$ which sends every generator to the single generator of $\mathbb{Z}$, we have that for every character $\chi \in \chara{G}$, $\Phi_{(\chi, \tilde{\alpha}, \sigma)}$ belongs to $\Iso (A\rtimes_{\tilde{\alpha}, r} G, (H\otimes \ell^2 (G))^{\oplus 2}, D_l)$.

\end{itemize}

\section{Postface}
As pointed out by a referee, it is not unnatural to equip quantum spaces with some sort of quantum symmetries.
Indeed, 
in the literature there is a notion of `quantum isometry group' for spectral triples (as well as for metric spaces) in the setting of (compact) quantum groups, see e.g. the monography \cite{GB}. Loosely speaking, these quantum isometry groups should be thought of as some sort of generalization of the $\Iso(A,H,D)$ considered in this paper, with groups replaced by quantum groups.
Even more interestingly, such quantum groups for the Christensen-Ivan spectral triples have been looked at in \cite{BGS} (cf. \cite{GB} Remark 5.2.2). Therein, the authors recognize the existence of some inductive limit structure on such quantum groups and present further considerations in the commutative case, notably for the algebra of continuous functions on the middle-third Cantor set. In general one would expect that the $\Iso$-groups associated to Christensen-Ivan spectral triples for $AF$-algebras considered in this work can be recovered as maximal classical subgroups (see e.g. \cite{ASS,KN}) of the corresponding quantum isometry groups.

We stress out the connection between Proposition \ref{prop3.1} and the the computation of the quantum isometry groups appearing in \cite{BGS}. In \cite{BGS} it is proved that the quantum isometry group $\mathsf{QISO}^+$ of a Christensen-Ivan spectral triple associated to the uniform measure on the Cantor set is an "inductive limit" of certain $C^*$-algebras $\mathsf{S}_n$ which are defined inductively by $\mathsf{S}_1 = C(\mathbb{Z}_2)$, $\mathsf{S}_{n+1} = (\mathsf{S}_n * \mathsf{S}_n) \oplus (\mathsf{S}_n * \mathsf{S}_n)$. Now, for every $n \in \mathbb{N}$, we let $\Iso_n (C(X), H_\mu, D)$ denote the (finite) group of automorphisms of $\mathcal{P}_n$-simple functions which respect the filtration of the algebra in the chosen direct limit decomposition. Then, as observed in the proof of Proposition \ref{prop3.1}, $\Iso_n (C(X), H_\mu , D) \simeq \ltimes_{k=1}^n S_2^{2^{k-1}}$ and $\Iso (C(X), H_\mu, D)$ is the projective limit of such groups. Using the fact that the abelianization of the free product of two groups is the direct product of the abelianizations, one can check that for every $n \in \mathbb{N}$ the maximal classical subgroup of $\mathsf{S}_n$ corresponds to $\Iso_n (C(X), H_\mu , D)$. The projective structure of $\Iso (C(X), H_\mu, D)$ is obtained by dualizing the inductive structure of $\mathsf{QISO}^+$. Hence $\Iso$ is the maximal classical subgroup of $\mathsf{QISO}^+$ in the sense of \cite{KN} (cf. \cite{D} Section 5).

\section{Acknowledgments}
Both authors are grateful to the anonymous referees for their valuable comments on the first draft of the paper and to A. Skalski for some clarification about the content of \cite{BGS}. They acknowledge the support of INdAM-GNAMPA. The first named author was supported by MIUR - Excellence Departments - grant: "$C^*$-algebras associated to $p$-adic groups, bi-exactness and topological dynamics", CUP: E83C18000100006, the grant Beyond Borders:  "Interaction of Operator Algebras with Quantum Physics and Noncommutative Structure", CUP: E84I19002200005 and the grant Beyond Borders: "A geometric approach to harmonic analysis and spectral theory on trees and graphs", CUP: E89C20000690005; he also acknowledges the University of Rome "Tor Vergata'' funding OAQM, CUP: E83C22001800005. R.C. is supported by the Sapienza Ricerca Scientifica 2019 Project "Algebre di operatori, analisi armonica, geometria noncommutativa ed applicazioni alla fisica quantistica, la combinatoria e la teoria dei numeri".

\bibliographystyle{mscplain}
 \bibliography{biblio}

\bigskip
{\parindent=0pt Addresses of the authors:\\

\smallskip Jacopo Bassi, Department of Mathematics, University of Tor Vergata, \\Via della Ricerca Scientifica 1, 00133 Roma, Italy. 
\\ E-mail: bssjcp01@uniroma2.it \\

\smallskip \noindent
Roberto Conti, Dipartimento SBAI,
Sapienza Universit\`a di Roma \\
Via A. Scarpa 16,
I-00161 Roma, Italy.
\\ E-mail: roberto.conti@sbai.uniroma1.it
\par}

\end{document}